\documentclass[psamsfonts 10pt]{amsart} 
\usepackage{amssymb, latexsym, eucal}
\usepackage[all]{xy}

\usepackage[utf8]{inputenc}

\usepackage{braket}
\usepackage{color}

\definecolor{darkgreen}{rgb}{0,.6,0}
\definecolor{darkred}{rgb}{.6,0,0}
\definecolor{darkblue}{rgb}{0,.2,.5}

\title[The graphs of join-semilattices and congruences of particle lattices]{The graphs of join-semilattices and the shape of congruence lattices of particle lattices}

\author[P. R\r{u}\v{z}i\v{c}ka]{Pavel R\r{u}\v{z}i\v{c}ka}

\address{
  Department of Algebra\\
	Faculty of mathematics and Physics\\
	Charles University in Prague\\
  Sokolovsk\'a 83\\
  186 75, Prague\\
  Czech republic
}
\email[P. R\r{u}\v{z}i\v{c}ka]{ruzicka@karlin.mff.cuni.cz}

\thanks{The author was partially supported by the Grant Agency of the Czech Republic under 
the grant no. GACR 14-15479S}
  


\date{\today}

\dedicatory{}

\numberwithin{equation}{section}

\theoremstyle{plain}                                     
\newtheorem{lemma}{Lemma}[section]
\newtheorem{theorem}[lemma]{Theorem}
\newtheorem{proposition}[lemma]{Proposition}
\newtheorem{corollary}[lemma]{Corollary}

\theoremstyle{remark}
\newtheorem{example}{Example}[section]

\theoremstyle{definition}

\newtheorem{claim}{Claim}

\newcommand{\qedc}{{\qed}~{\rm Claim~{\theclaim}.}}
\newenvironment{cprf} {\begin{proof}[Proof of Claim.]}
{\qedc\renewcommand{\qed}{}\end{proof}\setcounter{claim}{0}}

\newcommand{\setof}[2]{\ensuremath{\left\{ #1 \mid #2 \right\}}}                % Popis mnoziny
\newcommand{\fower}[1]{\ensuremath{\left[{#1}\right]^{<\omega}}}                % Mnozina konecnych podmnozin  
\newcommand{\idmap}[1]{\ensuremath{\boldsymbol{1}_{#1}}}												% Identity map

\newcommand{\poset}[1]{\ensuremath{{\boldsymbol{#1}}}}       % Znaceni usporadanych mnozin
\newcommand{\posi}{\poset{P}}
\newcommand{\posii}{\poset{P'}}
\newcommand{\posiii}{\poset{Q}}

\newcommand{\dset}[1]{{\downarrow\!(#1)}}
\newcommand{\uset}[1]{{\uparrow\!(#1)}}

\newcommand{\lattice}[1]{\ensuremath{{\boldsymbol{#1}}}}     % Typ znaceni svazu
\newcommand{\lati}{\lattice{L}}

\newcommand{\js}{\ensuremath{\langle 0,\vee \rangle}-semilattice}
\newcommand{\semilattice}[1]{\ensuremath{{\boldsymbol{#1}}}}     % Typ znaceni polo-svazu
\newcommand{\slati}{\semilattice{S}}

\newcommand{\con}[1]{\ensuremath{\mathrm{con}(#1)}}
\newcommand{\conS}[1]{\ensuremath{\mathrm{con}_{\vee}(#1)}}
\newcommand{\coni}{\ensuremath{\Theta}}
\newcommand{\conii}{\ensuremath{\Theta'}}

\newcommand{\econ}[1]{\ensuremath{\equiv_{#1}}}
\newcommand{\ncon}[1]{\ensuremath{\not \equiv_{#1}}}
\newcommand{\ecoi}{\econ{\coni}}

\newcommand{\ncoi}{\ncon{\coni}}

\newcommand{\lecon}[1]{\ensuremath{\leq_{#1}}}

\newcommand{\lecoi}{\ensuremath{\lecon{\coni}}}

\newcommand{\Graph}[1]{\ensuremath{{\boldsymbol{#1}}}}

\newcommand{\graphi}{\ensuremath{\Graph{G}_\lati}}

\newcommand{\graphsi}{\ensuremath{\Graph{G}_\slati}}

\newcommand{\joinirreducible}[1]{\ensuremath{J_{#1}}}

\newcommand{\jiri}{\joinirreducible{\lati}}

\newcommand{\jirsi}{\joinirreducible{\slati}}

\newcommand{\joinprime}[1]{\ensuremath{P_{#1}}}

\newcommand{\jpsi}{\joinprime{\slati}}

\newcommand{\edges}[1]{\ensuremath{E_{#1}}}
\newcommand{\edgi}{\edges{\lati}}

\newcommand{\edgsi}{\edges{\slati}}

\newcommand{\jirci}{\joinirreducible{\coni}}

\newcommand{\Joinirreducible}[1]{\ensuremath{J_{#1}^{\vee}}}

\newcommand{\Jirci}{\Joinirreducible{\coni}}
\newcommand{\Jircii}{\Joinirreducible{\conii}}

\newcommand{\Ideallattice}[1]{\ensuremath{\mathrm{Id}(#1)}}

\newcommand{\idsi}{\Ideallattice{\slati}}

\newcommand{\oid}[1]{\ensuremath{\mathcal{O}({#1})}}
\newcommand{\oideal}{\ensuremath{O}}

\newcommand{\hrd}[1]{\ensuremath{\mathrm{hrd}({#1})}}
\newcommand{\chrd}[1]{\ensuremath{\overline{\mathrm{hrd}}({#1})}}

\newcommand{\hint}[1]{\ensuremath{\partial(#1)}}

\newcommand{\hseti}{\ensuremath{H}}
\newcommand{\hsetii}{\ensuremath{H'}}

\newcommand{\conof}[1]{\ensuremath{\mathrm{\coni_{{#1}}}}}
\newcommand{\conH}{\conof{\hseti}}
\newcommand{\conHi}{\conof{\hsetii}}

\newcommand{\topol}[1]{\ensuremath{\mathcal{T}_{#1}}}
\newcommand{\topL}{\topol{\lati}}

\newcommand{\nbd}[2]{U({#1},{#2})}
\newcommand{\nbds}[1]{\mathcal{B}({#1})}

\begin{document}

\begin{abstract} We attach to each \js\ \slati\ a graph \graphsi\ whose vertices are join-irreducible elements of \slati\ and whose edges correspond to the reflexive dependency relation. We study properties of the graph \graphsi\ both when \slati\ is a join-semilattice and when it is a lattice. We call a \js\ \slati\ \emph{particle} provided that the set of its join-irreducible elements join-generates \slati\ and it satisfies DCC. We prove that the congruence lattice of a particle lattice is anti-isomorphic to the lattice of hereditary subsets of the corresponding graph that are closed in a certain zero-dimensional topology. Thus we extend the result known for principally chain finite lattices.       
\end{abstract}

\maketitle

\section{Introduction}\label{Introduction} % Section 1

The structure of congruences of  a finite lattice can be understood via the study of covers of their join-irreducible elements. The main tool to do so is the \emph{dependency relation} on the set of join-irreducible elements of the lattice (cf. \cite[page 39]{FJN95} or \cite[page 113]{Nat}). The idea of the use of the dependency relation goes back to \cite{PT75}, its original definition is due to A. Day \cite{Day79}. The dependency relation found a wide range of applications, aside of studying  congruences of finite \cite{FJN95}, resp. principally chain finite \cite{Nat}, lattices, let us mention characterization of finite lower bounded lattices. In particular, a finite lattice is lower bounded if and only if it contains no $D$-cycle \cite[Corollary 2.39]{FJN95}. 

A full description of a congruence lattice of a lattice via the dependency relation is established for finite lattice in \cite{FJN95} and generalized to principally chain finite lattices in \cite{Nat}. In this paper we extend these results to a wider class of \emph{particle} lattices, i.e., the lattices whose join-irreducible elements join-generate the lattice and they satisfy DCC (the decreasing chain condition). We also study properties of the dependency relation for \js{s}. In the \js\ case one cannot expect such a nice connection with the structure of the congruence lattice of the \js\ as in the case of lattices. Indeed, even when a \js\ is a small finite lattice, the lattice of its \js\ congruences may be much richer than the lattice of its lattice congruences (cf. Example~\ref{Ex:M3}). 

Let us sum up the content of the paper. Firstly we study consequences of DCC in posets. In particular we show that for posets satisfying DCC, saturated families of finite subsets of the posets satisfy certain minimality properties. This result is based on the fact that if a~poset satisfies DCC, the set of its anti-chains quasi-ordered by the \emph{join-refining relation} $\ll$ (see \cite[page 30]{FJN95}) satisfies DCC as well (cf. \cite[Exercise 10.4]{Nat}). In connection with these finiteness properties we study join-covers in join-semilattices and lattices. 

We define a graph of a \js. The set of vertices of the graph is the set of all join-irreducible elements of the \js\ and the edges correspond to the \emph{reflexive dependency relation} $\underline{D}$ defined in \cite[page 113]{Nat}. We prove that the graph has no edges other than loops for distributive semilattices and that it is symmetric for modular or relatively complement semilattices. This is a mild generalization of the corresponding results known for lattices (cf. \cite[Theorem 10.9]{Nat}).

We study how the congruence lattices of \js{s}, resp. lattices are related to the lattices of hereditary subsets of the corresponding graphs. We show that there is a Galois connection between the congruence lattice of a \js\ and the lattice of hereditary subsets of its graph. This connection is proved to be particularly nice for particle lattices. We define a zero-dimensional topology on the set of join-irreducible elements of a lattice, and we show that the Galois connection induces an anti-isomorphism between the congruence lattice of a particle lattice and the lattice of closed hereditary subsets of its graph. We apply this result to characterize congruence lattices of atomistic lattices.     

\section{Basic Concepts}\label{Basic Concepts} % Secton 2

Given a set $X$, we denote by $\fower{X}$ the set of all finite subsets of $X$. We denote by $\idmap{X}$ the identity map on the set $X$.

\subsection{Posets}\label{Posets} By a \emph{poset} we mean a partially ordered set. Given posets \posi\ and \posii, a map $f \colon \posi \to \posii$ is said to be \emph{monotone}, resp. \emph{antitone}, provided that $p \le q$ implies $f(p) \le f(q)$, resp. $p \le q$ implies $f(q) \le f(p)$, for all $p,q \in \posi$. An antitone bijection will be called \emph{anti-isomorphism}.

We say that a poset \posi\ satisfies DCC (\emph{the descending chain condition}) provided that there is no infinite decreasing sequence in \posi, equivalently, provided that each non-empty subset of \posi\ has a minimal element. Dually, we say that \posi\ satisfies ACC (\emph{the ascending chain condition}), if it does not contain an infinite increasing sequence. A subset $A \subseteq \posi$ is called an \emph{anti-chain} provided that the elements of $A$ are pairwise incomparable. 

A subset \oideal\ of a poset $\posi$ is said to be an \emph{order ideal} provided that $x \le y \in \oideal$ implies that $x \in \oideal$, for all $x,y \in \posi$. All order ideals of a poset form a sublattice of the lattice of all subsets of the poset. We denote the lattice of all order ideals of the poset \posi\ by \oid{\posi}\ (cf. Subsection~\ref{Strongly distributive lattices}).  

Given $X \subseteq \posi$, we set
\begin{equation*}
\begin{aligned}
\uset{X} &:= \setof{p \in \posi}{x \le p\ \text{for some}\ x \in X},\ \text{and dually}\\
\dset{X} &:= \setof{p \in \posi}{p \le x\ \text{for some}\ x \in X}.
\end{aligned}
\end{equation*} 
For a singleton set $X = \{x\}$, we abbreviate the notation writing $\dset{x}$ and $\uset{x}$. 

A \emph{tree} is a poset \poset{T} such that $\uset{x}$ is well-ordered for each $x \in \poset{T}$. The \emph{order type} $o(x)$ of an element $x \in \poset{T}$ is the order type of $\uset{x}$. For an ordinal $\alpha$, the \emph{$\alpha$th-level} of $\poset{T}$ is the set $\poset{T}_\alpha = \setof{x \in T}{o(x) = \alpha}$. A \emph{branch} of a tree is its maximal well-ordered subset. 

A \emph{quasi-order} on a set $Q$ is a binary relation, say $\ll$, on the set $Q$ such that $\ll$ is reflexive and transitive. Given $p,q \in Q$, we define $p \equiv q$ if both  $p \ll q$ and $q \ll p$. It is straightforward to verify that $\equiv$ is an equivalence relation on the set $Q$; we will call the relation $\equiv$ the \emph{equivalence induced by the quasi-order} $\ll$. For each $q \in Q$, we denote by $\overline{q}$ the block of the equivalence $\equiv$ containing $q$ and we set $\posiii = \setof{\overline q }{q \in Q}$. It is easy to see that the binary relation $\leq$ defined by $\overline p \le \overline q$ if and only if $p \ll q$ is well defined (i.e, it does not depend on the choice of the representatives of the blocks) partial order on the set \posiii. The poset \posiii\ will be called \emph{the maximal antisymmetric quotient} of $Q$.      

 Let $\posi, \posiii$ be posets. By a \emph{Galois connection} (between the posets \posi\ and \posiii) we mean a pair of antitone maps, $F \colon \posi \to \posiii$ and $G \colon \posiii \to \posi$, such that 
\begin{equation}\label{Eq:Galois1} 
p \le G(q)\ \text{if and only if}\ q \le F(p),\ \text{for all}\ p \in \posi, q \in \posiii.
\end{equation}    
Property \eqref{Eq:Galois1} is easily seen to be equivalent to 
\begin{equation}\label{Eq:Galois2} 
p \le GF(p)\ \text{for all}\ p \in \posi\ \text{and}\ q \le FG(q)\ \text{for all}\ q \in \posiii.
\end{equation} 

\subsection{Join-semilattices and lattices}\label{Join-semilattices} Let \slati\ be a join-semilattice. The least element of \slati, if it exists, will be denoted by $0$ and called \emph{zero} of the join-semilattice \slati. We will refer to join-semilattices with zero as \js. 

Let \slati\ be a join-semilattice (resp. a lattice). We say that $X \subseteq \slati$ \emph{join-generates} \slati\ provided that each element of \slati\ is the join of a finite subset of $X$. 

Let \slati\ be a \js. We say that $u \in \slati$ is \emph{join-irreducible} provided that $u = \bigvee F$ implies that $u \in F$, for all finite subsets $F$ of the \js\ \slati. Observe that since $0 = \bigvee \emptyset$, a join-irreducible element is necessarily non-zero. We denote by \jirsi\ the set of all join-irreducible elements of the \js\ \slati. 

An element $u$ of \js\ \slati\ is \emph{join-prime} provided that $u \le \bigvee F$ implies $u \le x$ for some $x \in F$, for all $F \in \fower{\slati}$. Clearly every join-prime element is join-irreducible, while join-irreducible elements may not be join-prime in general. We denote by $\jpsi$ the set of all join-prime elements of \slati. A non-zero element $u$ of \js\ \slati\ is \emph{completely join-prime} provided that $u \le \bigvee X$ implies $u \le x$ for some $x \in X$, for all (not necessarily finite) subsets $X$ of \slati. 

Let \slati\ be a \js, let $a,b \in \slati$. We say that $b$ \emph{covers} $a$, and we write $a \prec b$, if $a < b$ and $a \le x \le b$ implies $x \in \{a,b\}$ for all $x \in \slati$. An element $u$ of a \js\ \slati\ is an \emph{atom} provided that $0 \prec u$. A \js\ \slati\ is called \emph{atomistic} provided that each element of \slati\ is the join of a set of atoms (cf. \cite[p. 234]{Gra98}).  

An element $a$ of a complete lattice \lati\ is \emph{compact} provided that for all $X \subseteq \lati$, the inequality $a \le \bigvee X$ implies that $a \le \bigvee F$ for some finite $F \subseteq X$. An \emph{algebraic lattice} is a complete lattice whose each element is a (possibly infinite) join of compact elements.  

An \emph{ideal} of a join-semilattice \slati\ is its subset, say $I$, such that $x \vee y \in I$ if and only if both $x,y \in I$, i.e., an ideal of the semilattice \slati\ is its order ideal closed under finite joins. We denote by \idsi\ the lattice (necessarily algebraic) of all ideals of \slati.  

Given a join-semilattice \slati, we denote by \conS{\slati}\ the lattice of all congruences of \slati, i.e., equivalence relations \coni\ on \slati\ such that for all $a,b,c \in \slati$, $a \ecoi b$ implies that $a \vee c \ecoi b \vee c$. Give a lattice \lati, we denote by \con{\lati}\ the lattice of all congruence of \lati, i.e., equivalence relations on \lati\ respecting both the join and the meet. 

\subsection{Strongly distributive lattices}\label{Strongly distributive lattices} We say that a lattice \lati\ is \emph{strongly distributive} provided that it is isomorphic to the lattice $\oid{\posi}$ for some poset \posi. A strongly distributive lattice is necessarily algebraic and distributive, however not every algebraic distributive lattice is strongly distributive. Combining \cite[Lemma 10.6]{Nat} and \cite[Exercise 10.7]{Nat} we get that 
\begin{lemma}\label{L:StronglyDistributive} The following are equivalent for a distributive algebraic lattice \lati:
\begin{enumerate}
\item \lati\ is isomorphic to the lattice of order ideals of a poset.
\item Every element of \lati\ is a join of completely join-prime elements. 
\item Every compact element of \lati\ is a join of (finitely many) join irreducible compact elements.
\item The lattice \lati\ is dually algebraic.   
\end{enumerate}
\end{lemma}  

\subsection{Graphs}\label{Graphs} A \emph{graph} is a pair $\Graph{G} = (\joinirreducible{},\edges{})$ where $\joinirreducible{}$ is a set (of \emph{vertices}) and $\edges{} \subseteq \joinirreducible{} \times \joinirreducible{}$ (is a set of \emph{edges}). Given $u,v \in \joinirreducible{}$, we will write $u \to v$ to denote that $(u,v) \in \edges{}$. We say that $\hseti \subseteq \joinirreducible{}$ is \emph{hereditary} provided that if $u \in \hseti$ and $u \to v$, then $v \in \hseti$, i.e., a hereditary subset \hseti\ contains with each vertex $u \in \hseti$ all vertices reachable by an oriented path starting at $u$. 
\par
We denote by \hrd{\Graph{G}}\ the lattice of all hereditary subsets of $\joinirreducible{}$. Given $X \subseteq \joinirreducible{}$ we denote by \hint{X}\ the largest hereditary subset of $X$; equivalently, the union of all hereditary subsets of $X$. A subset $Y \subseteq \joinirreducible{}$ will be called \emph{co-hereditary} provided that its complement $\joinirreducible{} \setminus Y$ is hereditary.   

\section{Posets Satisfying DCC and minimal covers}

Let \posi\ be a poset, let $X$, $Y$ be subsets of \posi. We say that $X$ \emph{join-refines} $Y$, which we denote by $X \ll Y$, provided that $X \subseteq \dset{Y}$ \cite[p. 15]{FJN95}. By \cite[Lemma 1.15]{FJN95}, the relation $\ll$ forms a quasi-order on $\fower{\slati}$ and for each $X \in \fower{\posi}$, there is a unique anti-chain $A \in \fower{\posi}$ such that $A \equiv X$ (where $\equiv$ is the equivalence induced by the quasi-order $\ll$). Furthermore, if $A$ is an anti-chain with $A \equiv X$, then $A \subseteq X$.       

Let $\mathcal{C}$ be a set of finite subsets of \posi. We call the set $\mathcal{C}$ \emph{saturated} provided that for each non-empty $X \in \mathcal{C}$ there exists an anti-chain $A \in \mathcal{C}$ such that $A \subseteq X$. We say that $X \in \mathcal{C}$ is \emph{$\mathcal{C}$-minimal} provided that $Y \ll X$ implies $X \subseteq Y$, for all $Y \in \mathcal{C}$. 

\begin{lemma}\label{L:DCC} A poset $\posi = (P,\le)$ satisfies DCC if and only if each non-empty saturated $\mathcal{C} \subseteq \fower{\posi}$ contains a $ \mathcal{C}$-minimal element. 
\end{lemma}
     
\begin{proof} $(\Leftarrow)$ Suppose that the poset \posi\ does not satisfy DCC. Then there is an infinite strictly decreasing sequence $b_0 > b_1 > \cdots$ in \posi. Set $B = \{b_0,b_1,\cdots\}$ and $\mathcal{C} = \fower{B}$. One easily checks that $\mathcal{C}$ is saturated without a $\mathcal{C}$-minimal element. 
\par
$(\Rightarrow)$ Suppose that \posi\ satisfies DCC and there is a non-empty saturated $\mathcal{C} \subseteq \fower{\posi}$ with no $\mathcal{C}$-minimal element. 
\begin{claim}\label{Claim:DCC} There is an infinite sequence of pairwise distinct anti-chains $A_0,A_1,A_2,\dots$ in $\mathcal{C}$ such that 
\begin{equation}\label{Eq:DCC}
A_0 \gg A_1 \gg A_2 \gg \cdots.
\end{equation}
\end{claim}
\begin{cprf} We construct the sequence inductively. Since $\mathcal{C}$ is nonempty and saturated, there is a nonempty anti-chain $A_0 \in \mathcal{C}$. Let $n$ be a positive integer and suppose that we have already constructed a sequence $A_0,A_1,\cdots,A_{n-1}$ of pairwise distinct anti-chains satisfying \eqref{Eq:DCC}. By the assumption, there is no $\mathcal{C}$-minimal element. Therefore there is $X_n$ with $X_n \ll A_{n-1}$ and $A_{n-1} \nsubseteq X_n$. Since $A_{n-1}$ is an anti-chain, we have that $A_{n-1} \not \ll X_n$, hence $A_i \not \ll X_n$, whence $A_i \nsubseteq X_n$ for all $i = 0,\dots,n-1$. The set $\mathcal{C}$ is saturated, therefore $X_n$ contains a nonempty anti-chain $A_n \in \mathcal{C}$. From $A_n \subseteq X_n \ll A_{n-1}$, we get that $A_n \ll A_{n-1}$. Since $A_n \subseteq X_n$, we have that $A_n \neq A_i$ for all $i = 0,\dots,n-1$.   
\end{cprf}

Put $Q = \bigcup_{n = 0}^\infty A_n$ and observe that the set $Q$ is infinite. Indeed, it has infinitely many pairwise distinct subsets $A_n$. Let $\prec_Q$ be the \emph{covering relation} corresponding to the restriction of the order $\le$ to $Q$, i.e., $x \prec_Q z$ for $x,z \in Q$ if $x < z$ and $x < y < z$ for no $y \in Q$. Let $T$ be a set of all finite subsets $\{x_0,x_1,\cdots,x_k\}$ of $Q$ such that $x_0 \prec_Q x_1 \prec_Q \cdots \prec_Q x_k$, ordered by inclusion. Clearly $\poset{T} = (T,\subseteq)$ is an infinite tree and $o(\{x_0,x_1,\dots,x_k\}) = k$ for all $\{x_0,x_1,\dots,x_k\} \in T$. Since the sets $A_n$ are anti-chains, $x \prec_Q y$ implies that there is a (necessarily unique) positive integer $n$ with $y \in A_{n-1} \setminus A_n$ and $x \in A_n$. It follows that $\{x \in Q \mid x \prec_Q y\} \subseteq A_n$, in particular, the set is finite. Observing that $\poset{T}_0 = A_0$, it follows by induction that each level of the tree $\poset{T}$ is finite. By K$\ddot{\mathrm{o}}$nig's lemma \cite{Kon27}, the tree \poset{T}\ contains an infinite branch. This branch corresponds to an infinite strictly decreasing chain in $\posi$, which contradicts the assumption that \posi\ satisfies DCC.           
\end{proof}

Let \slati\ be a \js, $I \subseteq \slati$, and $a \in \slati$. An \emph{$I$-cover} of $a$ is a finite $F \subseteq I$ such that $a \le \bigvee F$. Let $\mathcal{C}_a$ denote the set of all $I$-covers of $a$. By a \emph{minimal $I$-cover of} $a$ we mean a $\mathcal{C}_a$-minimal element. A \emph{minimal join-cover of} $a$ is a minimal \slati-cover of $a$.  

Observe that given a poset \posi\ and $\mathcal{C} \subseteq \fower{\posi}$, each $\mathcal{C}$-minimal element is an anti-chain. This is true because a  $\mathcal{C}$-minimal element cannot contain a proper subset in $\mathcal{C}$, therefore it is an inti-chain due to $\mathcal{C}$ being saturated. In particular, if $I$ is a subset of a \js\ \slati\ and $a \in \slati$, each $I$-minimal cover of $a$ is an anti-chain. 

\begin{lemma}\label{L:MinimalInJoinDense} Let \slati\ be a \js\ and let $I$ be its join-generating subset. Then all minimal $I$-covers are minimal covers.
\end{lemma}

\begin{proof} Let $a \in \slati$ and let $F$ be a minimal $I$-cover of $a$. Let $G$ be a finite subset of \slati\ such that $a \le \bigvee G$ and $G \ll F$. Since $I$ is join-generating, $G$ refines to an $I$-cover of $a$, say $H$. Since $F$ is $I$-minimal, it follows that $F \subseteq H$, and so $F \ll H$. By the transitivity of $\ll$ we get that $F \ll G$, and since $F$ is an anti-chain, we conclude that $F \subseteq G$. Thus $F$ is a minimal join-cover of $a$.     
\end{proof}

We say that a \js\ \slati\ has the \emph{weak minimal join-cover refinement property} provided that each join-cover of an element $a \in \slati$ refines to a minimal join-cover of $a$. This is the minimal join-cover refinement property \cite[p. 30]{FJN95} weakened by removing the assumption that each element of \slati\ has only finitely many minimal join-covers. 

\begin{lemma}\label{L:MinimalJoinCoverRefinementProperty} Let \slati\ be a \js, let $I$ be its join-generating subset. Assume that $I$, viewed as a poset with the ordering inherited from \slati, satisfies DCC. Then \slati\ satisfies the weak minimal join-cover refinement property.     
\end{lemma}

\begin{proof} Let $a \in \slati$ and let $F$ be a join cover of $a$. Let $\mathcal{C}_a$ be the set of all $I$-covers of $a$ refining $F$. Since $I$  join-generates $\slati$, the set $\mathcal{C}_a$ is nonempty. Applying Lemma~\ref{L:DCC} and the assumption that $I$ satisfies DCC, we infer that there is a $\mathcal{C}_a$-minimal element. It is clearly an $I$-minimal cover of $a$. By Lemma~\ref{L:MinimalInJoinDense} all minimal $I$-covers of $a$ are minimal join-covers of $a$, and so we found a minimal join-cover of $a$ refining $F$.       
\end{proof}

Let us call a \js\ \slati\ \emph{particle} provided that \jirsi\ is join-generating subset of \slati\ satisfying DCC. A lattice is \emph{particle} if its join-semilattice reduct is a particle \js. We get readily from Lemma~\ref{L:MinimalJoinCoverRefinementProperty} that  

\begin{corollary}\label{Cor:Particle} A particle \js\ satisfies the weak minimal join-cover refinement property.
\end{corollary}

\begin{lemma}\label{L:DCCParticle} A \js\ satisfying DCC is particle.
\end{lemma}

\begin{proof} It suffices to prove that if \slati\ is a \js\ with DCC then \jirsi\ join-generates \slati. It is an easy exercise (see \cite[Lemma 2.8]{Nat}). 
\end{proof}

Note that the converse does not hold in general, indeed, any atomistic \js\ is particle but not all atomistic \js\ satisfy DCC. Let us finish this section with a partial converse of Lemma~\ref{L:MinimalJoinCoverRefinementProperty}, particularly implying that if \jirsi\ join-generates the \js\ \slati, the minimal join-covers and minimal \jirsi-covers coincide.   

\begin{lemma}\label{L:MinimalCovers} In a \js\ \slati\ each minimal cover is a minimal \jirsi-cover. 
\end{lemma}  

\begin{proof} To prove that each minimal cover is a minimal \jirsi-cover, it suffices to show that elements of minimal covers are join-irreducible. So let $F$ be a minimal join-cover of $a \in \slati$ and let $u \in F$. Suppose that $u$ is not join-irreducible. Then there exist $x,y \in \slati$ both strictly smaller than $u$ with $u = x \vee y$. Then $G = (F \setminus \{u\}) \cup \{x,y\}$ is a join-cover of $a$ such that $G \ll F$ and $F \not\subseteq G$. This contradicts the minimality of $F$.
\end{proof}

\begin{corollary}\label{C:MinimalCovers} If a \js\ \slati\ has the weak minimal join-cover refinement property, then \jirsi\ join-generates \slati. 
\end{corollary}

\section{The Graph of a Join-Semilattice}\label{The Graph of a Join-Semilattice} 

We define a \emph{graph of the \js} \slati\ to be the graph $\graphsi = (\jirsi,\edgsi)$, where the set $\edgsi$ of its edges is defined as follows: given $u,v \in \jirsi$, then $u \to v$ if there is $x \in \slati$ such that $u \le x \vee v$ but $u \nleq x \vee y$ for all $y < v$. Note that for principally join-finite lattices the relation $\edgsi$ corresponds to the reflexive dependency relation denoted in \cite[p. 113]{Nat} as $\underline{D}$. The next lemma is the semilattice version of \cite[Lemma~2.31]{FJN95}. 

\begin{lemma}\label{L:GraphWJCRP} Let \slati\ be \js\ satisfying the weak minimal join-cover refinement property, let $u,v \in \jirsi$. Then $u \to v$ if and only if $v$ belongs to a minimal join-cover of $u$. 
\end{lemma}   

\begin{proof} $(\Leftarrow)$ Let $F$ be a minimal join-cover of $u$ containing $v$. Put $x = \bigvee (F \setminus \{v\})$. Then $u \le \bigvee F = x \vee v$. The minimality of $F$ implies that $u \nleq x \vee y$ for all $y < v$. $(\Rightarrow)$ Suppose that $u \to v$ with $u \le x \vee v$ and $u \nleq x \vee y$ for all $y < v$. Then $\{x,v\}$ is a join-cover of $u$ and since \slati\ satisfies the weak minimal join cover refinement property, $\{x,v\}$ refines to a minimal join cover $F$ of $u$. Put $y = \bigvee \setof{z \in F}{z \le v}$ (note that this set is non-empty since $u \nleq x$). From $u \le x \vee y$, we get that $y = v$. Since $v$ is join-irreducible, we conclude that $v \in F$.   
\end{proof}

\begin{lemma}\label{L:GraphAndJoinPrimes} Let \slati\ be a \js, let $p \in \jpsi$ and $v \in \jirsi$. Then $p \to v$ implies that $p = v$.
\end{lemma}

\begin{proof} Suppose that $p \to v$. By the definition there is $x \in \slati$ such that $p \le x \vee v$ and $p \nleq x \vee y$ for all $y < v$. Since $p$ is join-prime, either $p \le x$ or $p \le v$. The first is impossible since $p \nleq x \vee 0$, the latter implies that $p = v$. 
\end{proof} 

Recall that a join-semilattice is \emph{distributive} provided that $a \le b \vee c$ implies that $a = y \vee z$ for some $y,z \in \slati$ with $y \le b$ and $z \le c$ (see e.g. \cite[p. 131]{Gra98}). The terminology comes from the fact that a join-semilattice is distributive if and only if its ideals form a distributive lattice \cite[Lemma II.5.1]{Gra98}.

It is straightforward that in a distributive \js, join-irreducible elements are join-prime. In fact, if the set of join-irreducible elements of a \js\ is join-generating, then join-prime and join-irreducible elements coincide if and only if the join-semilattice is distributive. Applying Lemma~\ref{L:GraphAndJoinPrimes}, we get readily that  

\begin{proposition}\label{P:Distributive} Let \slati\ be a distributive \js. Then $u \to v$ implies $u = v$ for all $u,v \in \jirsi$, i.e, the graph $\graphsi$ has no edges distinct from loops.  
\end{proposition} 

Following \cite{GS62} we say that a join-semilattice \slati\ is \emph{modular} provided that $a \le b \le a \vee c$ implies that there is $x \le c$ in \slati\ such that $b = a \vee x$ (see \cite{Rho75} for alternative definitions of modularity of semilattices). Similarly as in the previous case, a join-semilattice is modular if and only if its ideal lattice is modular.

\begin{proposition}\label{P:Modular} The graph of a modular \js\ \slati\ is symmetric, i.e, $u \to v$ implies $v \to u$ for all $u,v \in \jirsi$.
\end{proposition}

\begin{proof} Let $u,v$ be join-irreducible elements of \slati\ such that $u \to v$. By the definition of edges, there is $x \in \slati$ such that $u \le x \vee v$ but $u \nleq x \vee y$ for all $y < v$. 
\begin{claim}\label{Claim:Modular1} $v \le x \vee u$.
\end{claim}
\begin{cprf} Since $x \le x \vee u \le x \vee v$, there is $y \le v$ with $x \vee u = x \vee y$,  by the modularity. It follows that $u \le x \vee y$, hence $y = v$. We conclude that $v \le x \vee u$.  
\end{cprf}
Let $z \le u$ be such that $v \le x \vee z$. Then $z \le u \le x \vee v = x \vee z$, hence, by modularity, there is $w \le x$  with $u = w \vee z$. Since $u$ is join-irreducible, either $u = w \le x$, which is not the case, or $u \le z$. The latter means $u = z$. Thus we have proved that $v \to u$. 
\end{proof}
 
There is an alternative way how to prove Proposition~\ref{P:Modular}. Each \js\ \slati\ embeds into the lattice \idsi\ via the correspondence $a \mapsto \dset{a}$, sending each element of \slati\ to the corresponding principal ideal. It is straightforward to observe that $a \in \slati$ is join-irreducible if and only if the principal ideal $\dset{a}$ is join-irreducible in \idsi. 
 
\begin{lemma}\label{L:FullSubgraph} Let \slati\ be a \js, let $u,v \in \jirsi$. Then 
\begin{equation*}
u \to v\ \text{if and only if}\ \dset{u} \to \dset{v}.
\end{equation*} 
\end{lemma} 
 
\begin{proof} $(\Rightarrow)$ Suppose that $u \to v$. By the definition there is $x \in \slati$ such that $u \le x \vee v$ but $u \nleq x \vee y$ for all $y < v$. It follows that $\dset{u} \subseteq \dset{x}\ \vee \dset{v}$ and let $I \subseteq \dset{v}$ be an ideal of \slati\ such that $\dset{u} \subseteq \dset{x} \vee I$. Then there is $z \in I$ with $u \le x \vee z$. Since $I \subseteq \dset{v}$, we have that $z \le v$, hence $z = v$. It follows that $I = \dset{v}$, and so we have proved that $\dset{u} \to \dset{v}$.

$(\Leftarrow)$ Suppose that $\dset{u} \to \dset{v}$. Then there is an ideal of \slati, say $I$, such that $\dset{u} \subseteq I \vee \dset{v}$ and $\dset{u} \nsubseteq I \vee J$ for every ideal $J \subsetneq \dset{v}$. The first inequality implies that there is $x \in I$ with $u \le x \vee v$. Suppose that there is $y < v$ with $u \le x \vee v$. Then $\dset{y} \subsetneq \dset{v}$ and $\dset{u} \subseteq I \vee \dset{y}$. This is not the case, and so $u \to v$.    
\end{proof}

It follows from Lemma~\ref{L:FullSubgraph} that Proposition~\ref{P:Modular} reduces to the case when \slati\ is a modular lattices. In this case we can argue as in \cite[Theorem 10.9]{Nat}. 

Finally, we say that a \js\ \slati\ is \emph{relatively complemented} if for all $x \le y \le z$ in \slati\ there is $c \in \slati$ such that the meet $y \wedge c$ exists, $x = y \wedge c$, and $z = y \vee c$; we view the \js\ \slati\ as a partial lattice, assuming that the meet is defined whenever it exists. Note that \slati\ is not necessarily a lattice: Take for example the~lattice of all subspaces of an infinitely dimensional vector space $V$, consider its proper infinitely-dimensional subspace, say $W$, and remove all infinitely dimensional subspaces of $W$. The result is a relatively complemented \js\ that is not a lattice. 

\begin{lemma}\label{L:RelativelyComplemented} Let \slati\ be a relatively complemented \js, let $u,v \in \jirsi$. Then 
\begin{equation*}
u \to v \implies v \to u. 
\end{equation*} 
\end{lemma} 

\begin{proof} We can argue as in the proof of \cite[Theorem 10.9]{Nat}, observing that all join-irreducible elements of \slati\ are atoms.
\end{proof} 

\section{Congruences and join-irreducible elements}\label{Congruences and join-irreducible elements} % Section 5

Let \slati\ be a \js, let \coni\ be its congruence, and let $a,b \in \slati$. We write $a \ecoi b$ when $(a,b) \in \coni$, and $a \lecoi b$ when $(a \vee b,b) \in \coni$. Observe that $a \lecoi b$ is equivalent to $a \vee b \ecoi b$. Let us state and verify simple properties of these relations.

\begin{lemma}\label{L:SimpleProperties} Let \slati\ be a \js, let \coni\ be its congruence. Then the following holds true: 
\begin{enumerate}
\item For all $a,b \in \slati$, $a \ecoi b$ if and only if both $a \lecoi b$ and $b \lecoi a$.
\item The binary relation $\lecoi$ is a quasi-order on \slati.  
\end{enumerate}
\end{lemma}
 
\begin{proof}  (1) It is clear from the definition of the congruence relation on the semilattice that $a \ecoi b$ implies both $a \ecoi a \vee b$ and $b \ecoi a \vee b$. On the other hand $a \lecoi b$ is equivalent to $b \ecoi a \vee b$ and $b \lecoi a$ is equivalent to $a \ecoi a \vee b$. We conclude that $a \ecoi a \vee b \ecoi b$. (2) Let $a \lecoi b$ and $b \lecoi c$. It follows that $a \vee b \ecoi b$ and $b \vee c \ecoi c$, hence
\begin{equation*}
a \vee c \ecoi a \vee (b \vee c) = (a \vee b) \vee c \ecoi b \vee c \ecoi c,
\end{equation*}
whence $a \lecoi c$. Thus we proved transitivity of \lecoi, its reflexivity is clear from the definition.
\end{proof} 
 
Let \slati\ be a \js\ and let \coni\ be a congruence on \slati. We put
\begin{equation}\label{Eq:DefOfJv}
\Jirci := \setof{u \in \jirsi}{u \lecoi x \implies u \leq x,\ \text{for all}\ x \in \slati}.
\end{equation} 

\begin{lemma}\label{L:EquivDownset} Let \slati\ be a \js, let $\coni \in \conS{\slati}$, and let $a,b \in \slati$. Then
\begin{equation}\label{Eq:EquivDownSet1}
a \lecoi b \implies \dset{a} \cap \Jirci \subseteq \dset{b} \cap \Jirci.
\end{equation} 
\end{lemma}

\begin{proof} Let $u \in \dset{a} \cap \Jirci$. Since $u \le a$ and $a \lecoi b$, by the assumption, we infer, applying Lemma~\ref{L:SimpleProperties}(2), that 
$u \lecoi b$. Since $u \in \Jirci$, we conclude that $u \le b$, and so $u \in \dset{b} \cap \Jirci$. 
\end{proof}

Combining Lemmas~\ref{L:SimpleProperties}(1) and~\ref{L:EquivDownset}, we conclude that given a \js\ \slati, a congruence relation $\coni \in \conS{\slati}$, and elements $a,b \in \slati$, the implication
\begin{equation}\label{Eq:EquivDownSet2}
a \ecoi b \implies \dset{a} \cap \Jirci = \dset{b} \cap \Jirci
\end{equation}
holds true.

Given a \js\ \slati\ and a congruence relation $\coni \in \conS{\slati}$, we set
\begin{equation}\label{Eq:jirci}
\jirci := \setof{u \in \jirsi}{x < u \implies x \ncoi u,\ \text{for all}\ x \in \slati}.
\end{equation}
%Note that if we denote by $\lcoi$ the relation on \slati\ given by $a \lcoi b$ provided that $a \lecoi b$ and $a \ncoi b$, for all $a,b \in \slati$, then \eqref{Eq:jirci} can be reformulated as 
%\begin{equation*}
%\jirci := \setof{u \in \jirsi}{x < u \implies x \lcoi u,\ \text{for all}\ x \in \slati}.
%\end{equation*}
 
\begin{lemma}\label{L:DownSetEquiv} Let \slati\ be a particle \js, let $\coni$ be a congruence on \slati, and let $a,b \in \slati$. Then
\begin{equation}\label{Eq:DownSetEquiv}
\dset{a} \cap \jirci \subseteq \dset{b} \cap \jirci \implies a \lecoi b.  
\end{equation}
\end{lemma}

\begin{proof} The statement is clear when $a = 0$. Suppose that $0 < a$ and $\dset{a} \cap \jirci \subseteq \dset{b} \cap \jirci$. Put
\begin{equation*}
\mathcal{C}_a := \setof{A \in \fower{\jirsi}}{a \ecoi \bigvee A}.
\end{equation*}
One readily sees that $\mathcal{C}_a$ is saturated. Since \slati\ is a particle \js, the poset $\jirsi$ join-generates \slati\ and it satisfies DCC. Since  $\jirsi$ is join-generating in \slati, the set $\mathcal{C}_a$ is non-empty. Since the poset \jirsi\ satisfies DCC, there is an $\mathcal{C}_a$-minimal element, say $F$, due to Lemma~\ref{L:DCC}. 
\begin{claim}\label{Claim:DownSetEquiv}
The inclusion $F \subseteq \jirci$ holds true. 
\end{claim}
\begin{cprf}
Suppose the contrary. Then there is $y \in F$ such that $x \ecoi y$ for some $x < y$. If $x = 0$, then $\bigvee F \ecoi \bigvee \left(F \setminus \{y\} \right )$, hence $F \setminus \{y\} \in \mathcal{C}_a$, which contradicts the $\mathcal{C}_a$-minimality of $F$. If $0 < x$, then there is a finite $X \subseteq \jirsi$ with $x = \bigvee X$, (recall that \jirsi\ join-generates \slati). Put $G = (X \cup F) \setminus \{x\}$. From $x = \bigvee X$ we infer that $\bigvee G =  \bigvee F \ecoi a$. It follows that $G \in \mathcal{C}_a$. Observing that $G \ll F$ and $F \nsubseteq G$, (since $x \in F \setminus G$), we get the contradiction with the $\mathcal{C}_a$-minimality of $F$.   
\end{cprf}
From Claim~\ref{Claim:DownSetEquiv} we conclude that $F \subseteq \dset{a} \cap \jirci \subseteq \dset{b} \cap \jirci \subseteq \dset{b}$. It follows that $b \vee a \ecoi b \vee \bigvee F = b$, hence $a \lecoi b$. 
\end{proof}

Comparing the definitions of the sets $\jirci$ and $\Jirci$, we easily observe that $\Jirci \subseteq \jirci$. Indeed, $\jirci$ corresponds to the set of all join-irreducible elements minimal in their block $\coni$-blocks, while $\Jirci$ is the set of all join-irreducible elements which are minimum elements of their $\coni$-blocks. Observing that minimal elements of blocks of lattice congruences are necessarily unique in the blocks, hence minimum, we conclude that $\jirci = \Jirci$ when \lati\ is a lattice and $\coni \in \con{\lati}$.

\begin{corollary}\label{Corollary:DownsetEquiv}
Let \lati\ be particle lattice, let $\coni \in \con{\lati}$. Then for all $x,y \in \lati$,
\begin{equation*}
\dset{x} \cap \jirci \subseteq \dset{y} \cap \jirci \iff x \lecoi y.
\end{equation*}
\end{corollary}  

\begin{proof} Apply Lemmas~\ref{L:EquivDownset} and~\ref{L:DownSetEquiv}.
\end{proof}

\begin{lemma}\label{L:Arrow} Let \slati\ be a \js, let $\coni \in \conS{\slati}$, and let $u,v \in \jirsi$. Then the implication
\begin{equation}\label{Eq:Arrow}
\left( u \in \Jirci\ \text{and}\ u \to v \right) \implies v \in \jirci
\end{equation} 
holds true.
\end{lemma}

\begin{proof} Suppose that that there are $u,v \in \jirsi$ such that $u \to v$, $u \in \Jirci$, and $v \notin \jirci$. Since $u \to v$, there is $x \in \slati$ with $u \le x \vee v$ and $u \nleq x \vee y$ for all $y < v$. Since $v \notin \jirci$, there is $y \in \slati$ such that $y < v$ and $y \ecoi v$. The latter gives that $u \le x \vee v \ecoi x \vee y$. From $u \in \Jirci$ and $u \lecoi x \vee y$ we obtain that $u \le x \vee y$. This is a contradiction.   
\end{proof}

\begin{figure}[uht]
\begin{equation*}
\begin{xy}
(30,5)*{\bullet}*!U(3){0}, 
(10,15)*{\bullet}*!R(3){u},
(30,15)*{\bullet}*!R(3){x},
(50,15)*{\bullet}*!L(3){v},
(25,25)*{\bullet}*!UR(3){w},
(40,25)*{\bullet}*!DL(3){a},
(30,40)*{\bullet}*!D(3){1},
(30,5);(10,15) **@{-},
(30,5);(30,15) **@{-}, 
(30,5);(50,15) **@{-},
(30,15);(25,25) **@{-},
(30,15);(40,25) **@{-},
(50,15);(40,25) **@{-},
(10,15);(30,40) **@{-},
(40,25);(30,40) **@{-},
(25,25);(30,40) **@{-}
\end{xy}
\end{equation*}
\caption{The \js\ \slati.}
\end{figure}

\begin{example}\label{Example:JD} Consider the \js\ \slati\ depicted in Figure~1. Let $\coni \in \conS{\slati}$ be the least congruence identifying elements $u$ and $x$. The congruence \coni\ has exactly one non-singular block, namely  $\{a,u,w,x,1\}$. One easily observes that $\jirci = \{u,x,v\}$ and $\Jirci = \{v\}$. Since $u \le v \vee w$ but $u \nleq w$, we have that $u \to w$, $u \in \jirci$ and $w \not\in \jirci$. Similarly, since $v \le u \vee x$ and $v \nleq x$, we have that $v \to u$, $v \in \Jirci$, and $u \not\in \Jirci$. Therefore, the implication \eqref{L:Arrow} cannot be strengthen by either assuming that $u \in \jirci$ or concluding that $v \in \Jirci$. 
\end{example}

Of course, the situation simplifies when \coni\ is a lattice congruence. In this case Lemma~\ref{L:Arrow} corresponds to one implication of \cite[Theorem 10.5]{Nat}, (see also \cite[Lemma 2.33]{FJN95}). 

\begin{corollary}\label{Corollary:Arrow} Let \lati\ be a lattice, let $u,v \in \jiri$. Then for all $\coni \in \con{\lati}$:
\begin{equation*}
\left( u \in \jirci\ \text{and}\ u \to v \right) \implies v \in \jirci.
\end{equation*} 
\end{corollary}

\begin{lemma}\label{L:Hereditary} Let \slati\ be \js\ satisfying the weak minimal join-cover refinement property. Let $\hseti$ be a hereditary subset of \jirsi. Let $\conH$ be a binary relation on \slati\ defined by 
\begin{equation}\label{Eq:DefHereditary} a \equiv_{\conH} b \iff \dset{a} \cap \hseti = \dset{b} \cap \hseti \qquad (\text{for all}\ a,b \in \slati).
\end{equation}
Then $\conH$ is a congruence of the \js\ \slati\ preserving all existing meets. In particular, if \slati\ is a lattice, then $\conH \in \con{\slati}$.  
\end{lemma}

\begin{proof} It is clear from the definition, that the binary relation $\conH$ is reflexive, transitive, and symmetric, thus $\conH$ is an equivalence relation on \slati. Let $a,b,c \in \slati$ and suppose that $a \equiv_{\conH} b$. This, by \eqref{Eq:DefHereditary}, means that $\dset{a} \cap \hseti = \dset{b} \cap \hseti$. Let $u \in \dset{a \vee c} \cap \hseti$. Since \jirsi\ join-generates \slati\ due to Corollary~\ref{C:MinimalCovers}, we can find finite subsets $A$ and $C$ of  \jirsi\ such that $a = \bigvee A$ and $c = \bigvee C$. Observe that $A \cup C$ is a join-cover of $u$ and since \slati\ satisfies the minimal join-cover refinement property, $A \cup C$ refines to a minimal join-cover of $u$, say $F$. From Lemma~\ref{L:MinimalCovers} we get that $F \subseteq \jirsi$ and by Lemma~\ref{L:GraphWJCRP} we have that $u \to v$ for every $v \in F$. Since \hseti\ is hereditary, we infer that $F \subseteq \hseti$. Since $F \ll A \cup C$, either $v \subseteq a$ or $v \subseteq c$ for every $v \in F$. Since $\dset{a} \cap \hseti = \dset{b} \cap \hseti$, we have that $v \subseteq a$ implies $v \subseteq b$, hence $v \le b \vee c$, for all $v \in F$. It follows that $u \le \bigvee F \le b \vee c$. We conclude that $u \in \dset{b \vee c} \cap \hseti$.

Verifying that $\conH$ preserves existing meets is straightforward, indeed, for all $c \in \slati$, 
\begin{equation*}
\dset{c \wedge a} \cap H = \dset{c} \cap \dset{a} \cap H = \dset{c} \cap \dset{b} \cap H = \dset{c \wedge b} \cap H.
\end{equation*}     
\end{proof}

\section{The Galois connection} %Section 6

In this final section we study the connections between the congruence lattices of \js{s} (resp. lattices) and the lattices of all hereditary subsets of their graphs. We define a topology on the set of join-irreducible elements of a lattice \lati, induced by the ordering of \jiri, and we prove that the congruence lattice of a particle lattice is anti-isomorphic to the lattice of all closed hereditary subsets of its graph. Thus we generalize \cite[Corollary of Theorem 10.5]{Nat}. We apply this result to characterize the congruence lattices of atomistic lattices.

\begin{lemma}\label{L:GC}
Let \slati\ be a \js. The pair  $(F,G)$ of maps defined as
\begin{equation}\label{Eq:MapFG}
\begin{aligned}
F \colon  \conS{\slati} &\to \hrd{\graphsi} \\ 
 \coni &\mapsto \hint{\Jirci}. 
\end{aligned}
%\end{equation} and
\qquad \text{and} \qquad
%\begin{equation}\label{Eq:MapG}
\begin{aligned}
G \colon \hrd{\graphsi} &\to \conS{\slati} \\ 
 \hseti & \mapsto \conH,
\end{aligned}
\end{equation}
forms a Galois connection.
\end{lemma}

\begin{proof} First, let us carry out an easy verification of the antitonity of the maps $F$ and $G$. If $\coni \subseteq \conii$ for some $\coni,\conii \in \conS{\slati}$, then $\Jirci \supseteq \Jircii$ readily from definition \eqref{Eq:DefOfJv}. Consequently, $\hint{\Jirci} \supseteq \hint{\Jircii}$. To establish the latter, let $\hseti \subseteq \hsetii$ be hereditary subsets of \jirsi. It follows from definition \eqref{Eq:DefHereditary} that
\begin{equation*}
a \equiv_{\conHi} b \implies \dset{a} \cap \hsetii = \dset{b} \cap \hsetii \implies \dset{a} \cap \hseti = \dset{b} \cap \hseti \implies a \equiv_{\conH} b   
\end{equation*}
for all $a,b \in \slati$. Thus $\conH \supseteq \conH'$.   
\par
It remains to prove that for all $\coni \in \conS{\slati}$ and all $H \in \hrd{\graphsi}$:
\begin{equation*}
\coni \subseteq G(\hseti) \iff \hseti \subseteq F(\coni). 
\end{equation*}
\par $(\Rightarrow)$ Suppose that $\coni \subseteq G(\hseti) = \conH$. We are to prove that $\hseti \subseteq F(\coni) = \hint{\Jirci}$. Since $\hseti$ is a hereditary subset of $\jirsi$, it suffice to verify that $\hseti \subseteq \Jirci$. So let $u \in \hseti$, and let $a \in \slati$ be such that $u \lecoi a$. 
The assumption $\coni \subseteq \conH$ gives $u \lecon{\conH} a$, hence $\dset{u} \cap H \subseteq \dset{a} \cap H$. Since $u \in H$, we infer that $u \in \dset{a}$, that is, $u \le a$. We conclude that $u \in \Jirci$.  
\par $(\Leftarrow)$ Suppose that $\hseti \subseteq F(\coni) = \hint{\Jirci}$. Applying Lemma~\ref{L:EquivDownset}, we have for all $a,b \in \slati$ that
\begin{equation*}
a \ecoi b \implies \dset{a} \cap \Jirci = \dset{b} \cap \Jirci \implies \dset{a} \cap H = \dset{b} \cap H \implies a \econ{\conH} b.
\end{equation*}   
It follows that $\coni \subseteq \conH = G(\hseti)$.
\end{proof}

Let \lati\ be a lattice and let $\coni \in \con{\lati}$ be its congruence. As noted above $\jirci = \Jirci$, and by Corollary~\ref{Corollary:Arrow}, $\Jirci = \hint{\Jirci}$. On the other hand $G(\hseti) \in \con{L}$ for every hereditary subset $H$ of \jiri\ due to Lemma~\ref{L:Hereditary}. It follows that the Galois connection $(F,G)$, defined by \eqref{Eq:MapFG} above (now with the \js\ \slati\ replaced by the lattice \lati), consists of the pair of maps 
\begin{equation}\label{Eq:MapsFG}
\begin{aligned}
F \colon  \con{\lati} &\to \hrd{\graphi}  \\
 \coni &\mapsto \jirci 
\end{aligned}
\qquad \text{and} \qquad
\begin{aligned}
G \colon \hrd{\graphi} &\to \con{\lati} \\
\hseti & \mapsto \conH.
\end{aligned}
\end{equation}

\begin{lemma}\label{L:GF} Let $(F,G)$ be the Galois connection defined by \eqref{Eq:MapsFG}. If \lati\ is a particle lattice, then $GF = \idmap{\con{\lati}}$.
\end{lemma}

\begin{proof} By the definition, $GF(\coni) = \conof{\jirci}$ for all $\coni \in \con{\lati}$. Recalling definition \eqref{Eq:DefHereditary} and applying Corollary~\ref{Corollary:DownsetEquiv} 
we obtain the sequence of equivalences
\begin{equation*}
a \econ{\conof{\jirci}} b \iff \dset{a} \cap \jirci = \dset{b} \cap \jirci \iff a \ecoi b,\quad\text{for all } a,b \in \lati.
\end{equation*}
It follows that $\coni = \conof{\jirci}$, which was to prove.
\end{proof}

\begin{example}\label{Ex:M3} Even simple examples show that Lemma~\ref{L:GC} cannot be similarly extended when we consider the Galois connection between the congruence lattice of a particle (even finite) \js\ and the lattice of hereditary subsets of its graph. As an example consider the \js\ reduct of five-element modular non-distributive lattice \lattice{M_3}. The \js\  \lattice{M_3} has a graph with no proper non-trivial hereditary subset but its join-preserving congruences form the twelve element lattice on Figure~2.   

\begin{figure}[hbt]
\begin{equation*}
\begin{xy}
(20,-5)*{\lattice{M_3}},
(65,-5)*{\Graph{G}_{\lattice{M_3}}},
(110,-5)*{\conS{\lattice{M_3}}},
(20,5)*{\bullet}*!U(3){0}, 
(5,15)*{\bullet}*!UR(2){u},
(20,15)*{\bullet}*!R(3){v},
(35,15)*{\bullet}*!UL(2){w},
(20,25)*{\bullet}*!D(3){1},
(50,15)*{\bullet}*!UR(2){u},
(65,15)*{\bullet}*!U(3){v},
(80,15)*{\bullet}*!UL(2){w},
(110,5)*{\bullet}*!U(3){0},
(95,15)*{\bullet}, (110,15)*{\bullet}, (125,15)*{\bullet},
(95,25)*{\bullet}, (110,25)*{\bullet}, (125,25)*{\bullet},
(110,35)*{\bullet}, (105,35)*{\bullet}, (95,35)*{\bullet}, (125,35)*{\bullet},
(110,45)*{\bullet}*!D(3){1},
(20,5);(5,15) **@{-},
(20,5);(20,15) **@{-}, 
(20,5);(35,15) **@{-},
(20,25);(5,15) **@{-},
(20,25);(20,15) **@{-},
(20,25);(35,15) **@{-},
(110,5);(95,15) **@{-}, (110,5);(110,15) **@{-}, (110,5);(125,15) **@{-},
(95,15);(95,25) **@{-}, (95,15);(110,25) **@{-}, 
(110,15);(125,25) **@{-}, (110,15);(95,25) **@{-},
(125,15);(110,25) **@{-}, (125,15);(125,25) **@{-},
(110,35);(95,25) **@{-}, (110,35);(110,25) **@{-}, (110,35);(125,25) **@{-},
(95,35);(95,25) **@{-}, (105,35);(110,25) **@{-}, (125,35);(125,25) **@{-},
(95,35);(110,45) **@{-}, (105,35);(110,45) **@{-}, (110,35);(110,45) **@{-}, (125,35);(110,45) **@{-}, 
\ar(52,15);(62,15),\ar(62,16);(52,16),
\ar(67,15);(78,15),\ar(78,16);(67,16),
\ar@`{(65,25)} (52,17);(78,17),
\ar@`{(65,6)} (78,14);(52,14)
\end{xy}
\end{equation*}
\caption{The \js\ \lattice{M_3}.}
\end{figure}
\end{example} 

Let \lati\ be a lattice, let $x,v \in \lati$ be such that $v$ is join irreducible and $x < v$. For each such a pair of elements put
\begin{equation}\label{Eq:nbd}
\nbd{v}{x} = \setof{u \in \jiri}{u \nleq x \text{ and }u \le v} = \left( \dset{v} \setminus \dset{x} \right) \cap \jiri.   
\end{equation}
For each $v \in \jiri$ set 
\begin{equation*}
\nbds{v} = \setof{\nbd{v}{x}}{x < v\ \text{in}\ \lati}.
\end{equation*}

\begin{lemma}\label{L:Topology} Let \lati\ be a lattice. The collection $\{ \nbds{v} \}_{v \in \jiri}$ is a neighborhood system of a topological space.
\end{lemma}

\begin{proof} We shall verify that the collection $\{ \nbds{v} \}_{v \in \jirsi}$ satisfies properties (BP1-BP3) of \cite[page 12]{Eng89}. Let $v \in \jirsi$. Then $\nbds{v}$ is non-empty as $\nbd{v}{0} \in \nbds{v}$ and clearly $v \in U$ for each $U \in \nbds{v}$; thus (BP1) holds true. If $u \in \nbd{v}{x}$ for some $x < v$ in \slati\ and some $u \in \jirsi$, then $u \nleq x$ due to definition \eqref{Eq:nbd}, hence $u \wedge x < u$ and we have that $\nbd{u}{u \wedge x} \in \nbds{u}$. Since $u \le v$, we conclude that $\nbd{u}{u \wedge x} \subseteq \nbd{v}{x}$. This settles property (BP2). Finally, let $u \in \nbd{v_1}{x_1} \cap \nbd{v_2}{x_2}$, for some $u,v_i \in \jiri$ and $x_i < v_i$ in \lati, $i = 1,2$. Then, by definition \eqref{Eq:nbd}, we have that $u \nleq x_i$ and $u \le v_i$ for both $i = 1,2$. It follows that $u \wedge x_i < u$, $i = 1,2$, and since $u$ is join-irreducible, we conclude that $(u \wedge x_1) \vee (u \wedge x_2) < u$, and so $\nbd{u}{(u \wedge x_1) \vee (u \wedge x_2)} \in \nbds{u}$ is the desired neighborhood of $u$ with $\nbd{u}{(u \wedge x_1) \vee (u \wedge x_2)} \subseteq \nbd{v_1}{x_1} \cap \nbd{v_2}{x_2}$. This proves property (BP3).   
\end{proof}

Let \topL\ be the topology on the set \jiri\ generated by the neighborhood system $\bigcup_{v \in \jiri} \nbds{v}$ (cf. \cite[Proposition 1.2.3]{Eng89}). 

\begin{lemma}\label{P:Clopen} Let $\lati$ be a lattice, let $v \in \jiri$ and let $x < v$ in \lati. Then the neighborhood $\nbd{v}{x}$ is closed in the topology \topL. 
\end{lemma}

\begin{proof} Let $u \in \jirsi$ satisfy $u \notin \nbd{v}{x}$. If $u \nleq v$, then $u \wedge v < u$ and it is straightforward from \eqref{Eq:nbd} that $\nbd{u}{u \wedge v} \cap \nbd{v}{x} = \emptyset$. Suppose that $u \le v$. From $u \notin  \nbd{v}{x}$, we conclude that $u \le x$, hence $\nbd{u}{0} \cap \nbd{v}{x} = \emptyset$. It follows that $\jiri \setminus \nbd{v}{x}$ is open, and so $\nbd{v}{x}$ is closed.    
\end{proof}

Recall that a topological space is \emph{zero-dimensional} provided that it is $T_1$ and it has a basis consisting of \emph{clopen} sets, i.e., sets that are both closed and open, (see \cite[p. 360]{Eng89}). Let $\lati$ be a lattice. It follows readily from definition \eqref{Eq:nbd} that $\{v\} = \bigcap \nbds{v}$ for each $v \in \jiri$, thus the topology \topL\ is $T_1$. It follows from Lemma~\ref{P:Clopen} that the topology \topL\ has a basis of clopen sets, and so the topology is zero-dimensional.  

Given a lattice \lati, let $\chrd{\graphi}$ denote the lattice of all closed (w.r.t. the topology \topL) hereditary subsets of \jiri. 

\begin{lemma}\label{L:FG} Let \lati\ be a lattice and let $(F,G)$ be the Galois connection defined by \eqref{Eq:MapsFG}. Then the following holds true:
\begin{enumerate}
\item For every $\coni \in \con{\lati}$, $F(\coni) = \jirci \in \chrd{\graphi}$.
\item For every $\hseti \in \chrd{\graphi}$, $H = \joinirreducible{\conH}  = FG(H)$.    
\end{enumerate}
\end{lemma}

\begin{proof} $(1)$ Let $\coni \in \con{\lati}$ and let $v \in \jiri \setminus \jirci$. By definition \eqref{Eq:jirci}, there is $x \in \lati$ such that $x < v$ and $x \ecoi v$. We are going to show that the neighborhood $\nbd{v}{x}$ is included in $\jiri \setminus \jirci$. Let $u \in \nbd{v}{x}$, i.e., $u$ is a join-irreducible element of $\lati$ such that $u \le v$ and $u \nleq x$. Since $u \nleq x$, the inequality $x \wedge u < u$ holds true. From $x \ecoi v$ and $u \le v$, we infer that $x \wedge u \ecoi v \wedge u = u$. We conclude that $u \in \jiri \setminus \jirci$, and so the set $\jiri \setminus \jirci$ is open, hence its complement $\jirci$ is closed. Recall that \jirci\ is a hereditary subset of \jiri\ due to Corollary~\ref{Corollary:Arrow}.   
\par
(2) Let  $\hseti \in \chrd{\graphi}$. Since, by Lemma~\ref{L:GF}, the maps $(F,G)$ form a Galois connection, formula \eqref{Eq:Galois2} says that $H \subseteq \joinirreducible{\conH}$. In order to prove that $\joinirreducible{\conH} \subseteq H$, pick $u \in \joinirreducible{\conH}$ and let $x \in \lati$ be such that $x < u$. Then $x \ncon{\conH} u$, and so $\dset{x} \cap \hseti \subsetneq \dset{u} \cap \hseti$ due to \eqref{Eq:DefHereditary}, hence $\nbd{u}{x} = (\dset{u} \setminus \dset{x}) \cap \hseti \neq \emptyset$. It follows that $u$ belongs to the closure of \hseti. Since \hseti\ is supposed to be closed, we conclude that $u \in \hseti$. Thus we have proved the opposite inclusion $\joinirreducible{\conH} \subseteq H$. 
\end{proof}

\begin{theorem}\label{T:Main} Let \lati\ be a particle lattice and let $(F,G)$ be the Galois connection defined by \eqref{Eq:MapsFG}. Then the image of the map $F$ is $\chrd{\graphi}$, and the maps $F \colon \con{\lati} \to \chrd{\graphi}$ and $G' := G \restriction \chrd{\graphi} \colon \chrd{\graphi} \to \con{\lati}$ are mutually inverse lattice anti-isomorphisms. In particular, the lattice \con{\lati}\ is isomorphic to the lattice of all open co-hereditary subsets of \jiri.  \end{theorem}

\begin{proof} It follows from Lemma~\ref{L:FG} that the image of the map $F$ corresponds to \chrd{\graphi}\ and that $FG'  = \idmap{\chrd{\graphi}}$. If \lati\ is a particle lattice, we apply Lemmas~\ref{L:GF} and~\ref{L:FG}(1) to infer that $G'F = \idmap{\con{\lati}}$. We conclude that $F$ and $G'$ are mutually inverse lattice anti-isomorphisms. The last statement of the theorem easily follows.  
\end{proof}

\begin{proposition}\label{P:DownFinite} Let \lati\ be a lattice such that the set \jiri\ join-generating and $\dset{v}$ is finite for every $v \in \jiri$. Then the lattice \con{\lati}\ is strongly distributive.   
\end{proposition}

\begin{proof} The lattice \lati\ is clearly particle. Let $v \in \jiri$. Since $\dset{v}$ is finite, we infer that the singleton $\{v\} = \bigcap_{x < v} \nbd{v}{x}$ is open. It follows that the topology \topL\ is discrete. By Theorem~\ref{T:Main}, the lattice \con{\lati}\ is anti-isomorphic to the lattice of all hereditary subsets of the graph \graphi. Let $\ll$ be a transitive and reflexive closure of \edgi. Now it is straightforward to see that the lattice of all hereditary subsets of the graph \graphi\ is anti-isomorphic to the lattice of all order ideals of the maximal antisymmetric quotient of the quasi-ordered set 
$(\jiri,\ll)$. It follows that the lattice \con{\lati}\ is strongly distributive (cf. Lemma~\ref{L:StronglyDistributive}(4)). 
\end{proof}

It was proved by M. Tischendorf \cite{Tis92} that every finite lattice has a congruence preserving extension into an atomistic lattice. The construction of \cite{GS62}, adapted in \cite[Theorem 10.8]{Nat}, provides a representation of each algebraic strongly distributive lattice as \con{\lati} for a principally chain finite (i.e. no principal ideal of \lati\ contains an infinite chain \cite[p. 112]{Nat}) atomistic lattice. Applying Proposition~\ref{P:DownFinite} we prove that  

\begin{corollary}\label{C:Atomistic} Congruence lattice of atomistic lattice are exactly strongly distributive lattices. 
\end{corollary}  

In particular, Tischendorfs result cannot be extended beyond finite lattices.  

\bibliographystyle{amsplain}
\bibliography{myabbrev,bibwehrung,bibother}

\end{document}